\documentclass[a4paper,12pt]{article}
\usepackage{latexsym, amsmath, amsfonts, amscd, amssymb, verbatim, graphicx, graphics, fancyhdr}
\usepackage{color}
\usepackage{xcolor}

\def\tto{\;{\lower 1pt \hbox{$\rightarrow$}}\kern -10pt
\hbox{\raise 2pt \hbox{$\rightarrow$}}\;}

\def\proof{{\bf Proof.}}
\begin{document}
\pagestyle{myheadings}
\newtheorem{Theorem}{Theorem}[section]
\newtheorem{Proposition}[Theorem]{Proposition}
\newtheorem{Remark}[Theorem]{Remark}
\newtheorem{Lemma}[Theorem]{Lemma}
\newtheorem{Corollary}[Theorem]{Corollary}
\newtheorem{Definition}[Theorem]{Definition}
\newtheorem{Example}[Theorem]{Example}
\renewcommand{\theequation}{\thesection.\arabic{equation}}
\normalsize
\setcounter{equation}{0}

\title{\bf Approximate Solutions in Linear Fractional Vector Optimization}

\author{Nguyen Thi Thu Huong\footnote{Center for Applied Mathematics and Informatics, Institute of Information and Communication Technology, Le Quy Don Technical University, 236 Hoang Quoc Viet Road, Bac Tu Liem District, Hanoi, Vietnam; email: nguyenhuong2308.mta@gmail.com.}}\maketitle

\date{}

\medskip
\begin{quote}
\noindent {\bf Abstract.} 
This paper studies approximate solutions of a linear fractional vector optimization problem without requiring boundedness of the constraint set. We establish necessary and sufficient conditions for approximating weakly efficient points of such a problem via some properties of the objective function and a technical lemma related to the intersection of the topological closure of the cone generated by a subset of the Euclidean space and the interior of the negative orthant. As a consequence, we obtain necessary conditions and sufficient conditions for approximate efficient solutions to the considered problem. Applications of these results to linear vector optimization are considered.


\end{quote}

\section{Introduction}
\markboth{Approximate Efficiency in Linear Factional Vector Optimization}{\sc n. t. t. huong}
\setcounter{equation}{0}

\textit{Linear fractional vector optimization problems} (LFVOPs) are specific nonconvex vector optimization problems, which have been studied intensively due to their many noteworthy properties and theoretical importance (see \cite{B98,B01},\cite{C84}-\cite{CA83},\cite{HPY05a}-\cite{HHPY07},\cite{HYY_2020}-\cite{HY05}, \cite{Y2012,YY2011}, Chapter 8 of \cite{LTY05}, and the references therein).

 Connections of LFVOPs with monotone affine vector variational inequalities were firstly recognized by Yen and Phuong \cite{YP00}. Topological properties of the solution sets of LFVOPs and monotone affine vector variational inequalities have been studied by Choo and Atkins \cite{CA82,CA83}, Benoist~\cite{B98,B01}, Huy and Yen \cite{HY05}, Hoa et al.~\cite{HPY05a,HPY05b,HHPY07}, Huong et al.~\cite{HYY_2020,HYY_Optim2017} and other authors. By a fundamental theorem on stability of monotone affine
variational inequalities, Yen and Yao \cite{YY2011} derived several results on solution
stability and topological properties of the solution sets of LFVOPs. Then, Yen and Yang~\cite{Yen_Yang_JOTA2011} initiated a study on infinite-dimensional LFVOPs via affine variational inequalities on normed spaces. Numerical methods for solving LFVOPs can be found in Malivert \cite{M95} and Steuer~\cite{S86}. The interested readers are referred to the survey paper of Yen~\cite{Y2012} for more information about linear fractional vector optimization problems.

The notions of efficient solutions and weakly efficient solutions are crucial for analyzing vector optimization problems (see, e.g., \cite{Luc89,S86}). But, sometimes the corresponding solution sets are empty. In that case, one may wish  to find approximate solutions, which satisfy some requirements of the decision maker. In addition, note that some algorithms such as iterative algorithms and search algorithms often provide approximate solutions. Therefore, considering approximate solutions and studying their necessary and sufficient conditions is an important question from both theoretical and practical points of view. Based on the Kutateladze's concept~\cite{K79} of approximate points, Loridan \cite{Loridan84} introduced the notion of $\varepsilon$-efficient solutions to vector optimization and obtained some similar results as in~\cite{K79}. Then, several authors (see, e.g.~\cite{CK2016,H_2024,LiWang,Loridan84,ST14,White86} and the references therein) have established further results in this direction. 

Li and Wang~\cite{LiWang} proposed the concept of $\varepsilon$-proper efficiency in vector optimization and obtained several necessary and sufficient conditions for $\varepsilon$-proper efficiency via scalarization and an alternative theorem. Afterwards, Liu~\cite{Liu99} obtained some scalarization results for the kind of approximate properly efficient solutions for general vector optimization problems. Recently, for a linear fractional vector optimization problem with a bounded constraint set, Tuyen~\cite[Theorem~3.2]{T22} has shown that there is no difference between the $\varepsilon$-properly efficient solution set and the $\varepsilon$-efficient solution one. Later, in \cite[Theorem~3.3]{H_2024} we have presented necessary and sufficient conditions for a $\varepsilon$-proper efficient solution of a general vector optimization problem via Benson's approach~\cite{Be79}. Moreover, we have shown~\cite[Theorem~5.4]{H_2024} that for any linear vector optimization problem with a pointed polyhedral convex cone $K$, either the $e$-properly efficient solution set is empty or it coincides with the $e$-efficient solution set, where $e$ is any nonzero vector taken from the closed pointed ordering cone. 
 
It is well known that for a general vector optimization problem, the $\varepsilon$-properly efficient solution set is a subset of the $\varepsilon$-efficient solution set and the $\varepsilon$-efficient solution set is a subset of the $\varepsilon$-weakly efficient solution set.

In this paper, we establish necessary and sufficient conditions for approximating weakly efficient points of a linear fractional vector optimization problem without requiring the constraint set's boundedness. As a consequence, we obtain necessary conditions and sufficient conditions for approximate efficient solutions to the considered problem. These results can be applied to linear vector optimization problems.

Section~\ref{Sect_2} recalls some definitions and auxiliary results. Section~\ref{Sect_3} establishes the main results and presents an illustrative example which has no weakly efficient solution, while both approximate efficient solution set and approximate weakly efficient solution set are nonempty.

\section{Preliminaries}\label{Sect_2}
\markboth{Approximate Proper Efficiency in Vector Optimization}{\sc n. t. t. huong}
\setcounter{equation}{0}
The scalar product and the norm in the Euclidean space $\mathbb R^p$ are denoted, respectively, by $\langle\cdot,\cdot\rangle$ and $\|\cdot\|$. 
Vectors in $\mathbb R^p$ are represented by columns of real numbers. If $A$ is a matrix, then $A^T$ denotes the transposed matrix of $A$. Thus, one has $\langle x,y\rangle=x^Ty$  for any $x,y\in \mathbb R^p$. For $x=(x_1,\ldots,x_p)$ and $y=(y_1,\ldots,y_p)$ from $\mathbb R^p$, one writes  $x\leq y$ (resp., $x<y$) whenever $x_i\leq y_i$ (resp., $x_i< y_i$) for all $i=1,\ldots,p$. The nonnegative orthant in $\mathbb R^p$ and the set of positive integers are denoted respectively by $\mathbb R^p_+$ and $\mathbb N$. 

A nonempty set $K\subset\mathbb R^m$ is called a \textit{cone} if $tv\in K$ for all $v\in K$ and $t\geq 0$. One says that $K$ is \textit{pointed} if $K\cap (-K)=\{0\}$. The smallest cone containing a nonempty set $D\subset\mathbb R^m$, i.e., the cone generated by $D$, will be denoted by ${\rm cone}\,D$. The topological closure of $D$ is denoted by $\overline{D}$ and $\overline{\rm cone}\,D:=\overline{{\rm cone}\, D}$.

A nonzero vector $v\in\mathbb R^n$ (see \cite[p.~61]{Roc70}) is said to be a \textit{direction of recession} of a nonempty convex set $D\subset\mathbb R^n$ if $x+tv\in D$ for every $t\geq 0$ and every $x\in D.$  The set composed by $0\in\mathbb R^n$ and all the directions $v\in\mathbb R^n\setminus\{0\}$ satisfying the last condition, is called the \textit{recession cone} of $D$ and denoted by $0^+D.$  If $D$ is closed and convex, then  $$0^+D=\{v\in\mathbb R^n\;:\; \exists x\in\Omega\ \, {\rm s.t.}\ \, x+tv\in D\ \, {\rm for\ all}\ \, t> 0\}.$$

\medskip
Consider \textit{linear fractional functions} $f_i:\mathbb R^n\to\mathbb R,\ i=1,\dots,m$,
of the form
$$f_i(x)=\frac{a_i^Tx+\alpha_i}{b_i^Tx+\beta_i},$$ where $a_i\in\mathbb R^n, b_i\in\mathbb   R^n, \alpha_i\in\mathbb  R,$
and $\beta_i\in\mathbb R$. Let $K$ be a \textit{polyhedral convex set}, i.e., there exist $p\in \mathbb N$, a matrix $C=(c_{ij})\in {\mathbb R}^{p\times n},$ and a vector $d=(d_i)\in {\mathbb R}^p$ such that $K=\big\{x\in {\mathbb R}^n\, :\, Cx\leq d\big\}$. Our standing condition is that
$b_i^Tx+\beta_i>0$ for all $i\in I$ and $x\in K$, where $I:=\{1,\dots,m\}$. Put
$f(x)=(f_1(x),\dots,f_m(x))$ and let $$\Omega=\big\{x\in\mathbb R^n\,:\,b_i^Tx+\beta_i>0,\ \, \forall i\in I\big\}.$$ Clearly, $\Omega$ is open and convex, $K\subset\Omega$, and $f$ is continuously differentiable on $\Omega$. The {\it linear fractional vector optimization problem} (LFVOP) given by $f$, $K$, and the ordering cone $\mathbb R^m_+$, is formally written as

\smallskip
\hskip1cm (VP) \hskip 2cm ${\rm Minimize}\ \, f(x)\ \,$
subject\ to $\; x\in K$.

\begin{Definition} {\rm A point $x\in K$ is said to be an {\em efficient solution} (or a {\em Pareto solution}) of ${\rm (VP)}$ if $\big(f(K)-f(x)\big)\cap \big(-\mathbb R^m_+\setminus\{0\}\big)=\emptyset$. One calls $x\in K$ a {\em weakly efficient solution} (or a {\em weak Pareto solution}) of ${\rm (VP)}$ if $\big(f(K)-f(x)\big)\cap \big(-{\rm int}\,\mathbb R^m_+\big)=\emptyset$, where ${\rm int}\,\mathbb R^m_+=\big\{\xi=(\xi_1,\dots,\xi_m)\in\mathbb R^m\,:\,\xi_i>0\ {\rm for}\ i=1,\dots,m\big\}$ denotes the interior of~$\mathbb R^m_+$.}
\end{Definition}

The efficient solution set (resp., the weakly efficient solution set) of (VP) are denoted,  respectively, by $E$ and $E^w$. According to \cite{CA82,M95} (see also \cite[Theorem 8.1]{LTY05}), \textit{for any $x\in K$, one has $x\in E$ (resp., $x\in E^w$) if and only if there exists a multiplier $\xi=(\xi_1,\dots,\xi_m)\in {\rm int}\,\mathbb R^m_+$ (resp., $\xi=(\xi_1,\dots,\xi_m)\in\mathbb R^m_+\setminus\{0\}$) such that
	\begin{equation*}\label{Eq8.3}\Big\langle \sum_{i=1}^m\xi_i\left[\big(b_i^Tx
		+\beta_i\big)a_i-\big(a_i^Tx+\alpha_i)b_i\right],y-x\Big\rangle\geq
		0, \quad  \forall y\in K.\end{equation*}}

If $b_i=0$ and $\beta_i=1$ for all $i\in I$, then (VP) coincides with the classical \textit{linear vector optimization problem} (LVOP) (see Luc~\cite{Luc_Book2016} and the references therein). By the above optimality conditions, for any $x\in K$, one has $x\in E$ (resp., $x\in E^w$) if and only if there exists a multiplier $\xi=(\xi_1,\dots,\xi_m)\in {\rm int}\,\mathbb R^m_+$ (resp.,  a multiplier $\xi=(\xi_1,\dots,\xi_m)\in\mathbb R^m_+\setminus\{0\}$) such that
\begin{equation*}\label{Eq8.3(1)}\Big\langle \sum_{i=1}^m\xi_ia_i,y-x\Big\rangle\geq
	0, \quad  \forall y\in K.\end{equation*}

Let $\varepsilon=(\varepsilon_1,...,\varepsilon_m)$ be a vector in $\mathbb{R}^m_+$. Specializing the concept of $\varepsilon$-efficiency of general vector optimization problems in~\cite{LiWang,Liu99,Loridan84} to (VP) we have the next definition.

\begin{Definition}\label{appro_def_efficiency}{\rm A point $x\in X$ is said to be an {\it $\varepsilon$-efficient solution} of (VP) if there exists no $y\in K$ such that $f(y)\leq f(x)-\varepsilon$ and $f(y)\neq f(x)-\varepsilon$.} 
\end{Definition}

Slightly weakening the requirement of $\varepsilon$-efficiency, we get the following notion of $\varepsilon$-weak efficiency.

\begin{Definition}\label{appro_def_weak_eff}{\rm A point $x\in K$ is said to be an {\it $\varepsilon$-weakly efficient solution} of (VP) if there exists no $y\in K$ such that $f(y)<f(x)-\varepsilon$.} 
\end{Definition}

The $\varepsilon$-efficient solution set (resp., the $\varepsilon$-weakly efficient solution set) of (VP) are denoted,  respectively, by $E_{\varepsilon}$ and $E_{\varepsilon}^{w}$. 

\begin{Remark}\label{Remark 1} {\rm Let $\bar x \in K$. Then, $\bar x$ is an $\varepsilon$-efficient solution of (VP) if and only if $\big[f(K)-(f(\bar x)-\varepsilon)\big]\cap\big(-\mathbb R^m_+\setminus\{0\}\big)=\emptyset$. Similarly, $\bar x$ is an $\varepsilon$-weakly efficient solution of (VP) if and only if $\big[f(K)-(f(\bar x)-\varepsilon)\big]\cap \big(-{\rm int}\,\mathbb R^m_+\big)=\emptyset$. Clearly, $E_{\varepsilon}\subset E_{\varepsilon}^{w}$.}
\end{Remark}

When $\varepsilon=0$, the notion of $\varepsilon$-efficient solution (resp., the notion of $\varepsilon$-weakly efficient solution) reduce to the notion of efficient solution, respectively, the notion of weakly efficient solution, i.e., $E_0=E$ and $E_0^w=E^w$.

In the sequel, to establish verifiable necessary sufficient conditions for a point $\bar x\in K$ to belong to $E^w_\varepsilon$, we will need the two following lemmas. The first one is a fact related to the intersection of the topological closure of the cone generated by a subset of $\mathbb{R}^m$ and the interior of the negative orthant.

\begin{Lemma}\label{Lemma2} 
Let $\Omega$ be a nonempty subset of $\mathbb{R}^m$. Then the following two properties are equivalent: 
\begin{description}
\item{\rm (i)} $\Omega\cap \big(-{\rm int}\,\mathbb R^m_+\big) =\emptyset$;
\item{\rm (ii)} $\overline{\rm cone}\,\Omega\cap \big (-{\rm int}\,\mathbb R^m_+\big)=\emptyset$.
\end{description}
\end{Lemma}
\proof\ Since $\Omega\subset {\rm cone}\,\Omega\subset \overline{\rm cone}\,\Omega$ by definition,~(ii) implies~(i). To prove the reverse implication, we can argue by contradiction. Suppose that~(i) holds, but~(ii) is invalid Then, there exists a vector $\bar v=(\bar v_1,...,\bar v_m)\in-{\rm int}\,\mathbb R^m_+$ such that $\bar v=\displaystyle\lim_{k\rightarrow\infty}v^k$ with $v^k=t_kx^k$, for  $x^k\in\Omega$ and $t_k>0$ for all $k\in\mathbb N$. As $\displaystyle\lim_{k\rightarrow\infty}v^k_i=\bar v_i<0$ for every $i\in I$, there exists $\bar k\in\mathbb N$ such that $v^k_i<0$ for all $k>\bar k$ and $i\in I$. So, we get 
$$x^k\in \Omega\cap \left(-{\rm int}\,\mathbb R^m_+ \right)\quad \forall k\geq\bar k,$$ a contradiction to~(i).
$\hfill\Box$ 

\begin{Lemma}[{See, e.g., \cite[Lemma~8.1]{LTY05} and \cite{M95}}]\label{Lemma1_lff} 
	Let $\varphi(x)=\dfrac{a^Tx+\alpha}{b^Tx+\beta}$ be a linear fractional function defined by $a,b\in\mathbb R^n$ and $\alpha,\beta\in\mathbb R$. Suppose that  $b^Tx+\beta\neq 0$ for every $x\in K_0$, where $K_0\subset\mathbb R^n$ is an arbitrary convex set. Then, one has
	\begin{equation*}\label{derivative}
		\varphi(y)-\varphi(x)=\frac{b^Tx+\beta}{b^Ty+\beta}\, \langle\nabla \varphi(x),y-x\rangle,
	\end{equation*} for any $x, y\in K_0$, where $\nabla \varphi(x)$
	denotes the Fr\'echet derivative of $\varphi$ at $x$.
\end{Lemma}

\newpage
\section{Necessary and sufficient conditions for $\varepsilon$-efficiency}~\label{Sect_3}
\markboth{Approximate Proper Efficiency in Vector Optimization}{\sc n. t. t. huong}
\setcounter{equation}{0}

First, necessary and sufficient conditions for a feasible point to be an  $\varepsilon$-weakly efficient solution of~(VP), where $\varepsilon\in\mathbb{R}^m_+$, are provided by the following theorem.

\begin{Theorem}\label{Proposition-1}  Let $\varepsilon\in\mathbb{R}^m_+$ and $\bar x\in K$. Then, $\bar x$ is an $\varepsilon$-weakly efficient solution of~{\rm (VP)} if and only if there exists a vector $\lambda=(\lambda_1,\ldots,\lambda_m)\in\mathbb{R}^m_+\setminus \{0\}$ such that
	\begin{equation}\label{equ-2}
		\sum_{i=1}^{m}\lambda_i[(b_i^T\bar x+\beta_i)\langle\nabla f_i(\bar x), y-\bar x\rangle+ \varepsilon_i(b_i^Ty+\beta_i)]\geq 0
	\end{equation}
	for all $y\in K$.
\end{Theorem}
\proof\ 
	By definition, $\bar x$ is an $\varepsilon$-weakly efficient solution of (VP) if and only if the following system has no  $y\in K$ such that
	\begin{equation}\label{system-1}
		f_i(y)< f_i(\bar x)-\varepsilon_i, \ \ \forall i\in I.
	\end{equation}
	For every $i\in I$, by Lemma~\ref{Lemma1_lff} one has 
	$$f_i(y)-f_i(\bar x)=\frac{b_i^T\bar x+\beta_i}{b_i^T y+\beta_i}\langle \nabla f_i(\bar x), y-\bar x\rangle.$$ So, the system~\eqref{system-1} of inequalities can be rewritten as follows:
			$$\frac{b_i^T\bar x+\beta_i}{b_i^T y+\beta_i}\langle \nabla f_i(\bar x), y-\bar x\rangle+\varepsilon_i< 0,  \ \ \forall i\in I.$$
Since $b_i^Ty+\beta_i>0$ for all $i\in I$ and $y\in K$, the latter is equivalent to the condition
	\begin{equation*}\label{system-2}
(b_i^T\bar x+\beta_i)\langle \nabla f_i(\bar x), y-\bar x\rangle+\varepsilon_i(b_i^T y+\beta_i)< 0, \ \ \forall i\in I.  
	\end{equation*} 
	Let
	\begin{equation*} 
		A:=
		\begin{pmatrix}
			(b_1^T\bar x+\beta_1)\nabla f_1(\bar x)^T+\varepsilon_1b_1^T
			\\
			\vdots
			\\
			(b_m^T\bar x+\beta_m)\nabla f_m(\bar x)^T+\varepsilon_mb_m^T
		\end{pmatrix}, \ \
		b:=
		\begin{pmatrix}
			-(b_1^T\bar x+\beta_1)\langle \nabla f_1(\bar x), \bar x\rangle+\varepsilon_1\beta_1
			\\
			\vdots
			\\
			-(b_m^T\bar x+\beta_m)\langle \nabla f_m(\bar x), \bar x\rangle+\varepsilon_m\beta_m
		\end{pmatrix}.
	\end{equation*}
	Then, $\bar x$ is an $\varepsilon$-weakly efficient solution of~(VP) if and only if \textit{there is no $y\in K$ such that $Ay+b \in (-{\rm int}\,\mathbb R^m_+)$}. The last condition means that
	\begin{equation}\label{equa-1}
		D\cap(-{\rm int}\,\mathbb R^m_+)=\emptyset, 
	\end{equation} where $D:=\{Ay+b\,:\, y\in K\}.$ Put $P=\mathrm{cone}\,D$ and observe that~\eqref{equa-1} can be rewritten as $P\cap(-{\rm int}\,\mathbb R^m_+)=\emptyset$. So, by Lemma~\ref{Lemma2},~\eqref{equa-1} is equivalent to the condition 
\begin{equation}\label{equa-8}
	\overline P\cap(-{\rm int}\,\mathbb R^m_+)=\emptyset. 
\end{equation}  
Let $P^*:=\{z^*\in\mathbb{R}^m\,:\, \langle z^*, z\rangle\geq 0,\; \forall z\in P\}.$
We claim that~\eqref{equa-8} holds if and only if
	\begin{equation}\label{equ-9}
		P^*\cap \big(\mathbb{R}^m_+\setminus \{0\}\big)\neq\emptyset.  
	\end{equation}
	Indeed, arguing by contradiction, one can easily show that~\eqref{equ-9} implies~\eqref{equa-8}. We now assume that~\eqref{equa-8} holds. Then, applying the separation theorem \cite[Theorem~11.3]{Roc70} to the convex sets $\overline P$ and  $(-\mathbb R^m_+)$ yields a vector $\xi\in\mathbb{R}^m\setminus\{0\}$ such that  
	\begin{equation}\label{eq_10}
		\langle\xi, u\rangle\leq\langle\xi, v\rangle, \ \ \forall u\in -\mathbb{R}^m_+,\ \ \forall v\in {\overline P}.
	\end{equation} 
Fixing any  $v\in {\overline P}$, from~\eqref{eq_10} we can deduce that $\xi\in\mathbb{R}^m_+\setminus \{0\}$. Moreover, substituting $u=0$ to the inequality in~\eqref{eq_10} shows that $\xi\in P^*$. So,~\eqref{equ-9} holds. 
	
Clearly,~\eqref{equ-9} means there exists a vector $\lambda=(\lambda_1,\ldots,\lambda_m)\in\mathbb{R}^m_+\setminus \{0\}$ such that $\langle\lambda, v\rangle\geq 0$ for all $v\in\bar P$. Since the latter can be rewritten as  
	\begin{equation*}
		\langle\lambda, Ay+b\rangle\geq 0\quad \forall y\in K
	\end{equation*} and, by the constructions of $A$ and $b$,
$$\langle\lambda, Ay+b\rangle=\sum_{i=1}^{m}\lambda_i[(b_i^T\bar x+\beta_i)\langle\nabla f_i(\bar x), y-\bar x\rangle+ \varepsilon_i(b_i^Ty+\beta_i)],$$ we have thus proved that $\bar x$ is an $\varepsilon$-weakly efficient solution of~{\rm (VP)} if and only if there exists a vector $\lambda=(\lambda_1,\ldots,\lambda_m)\in\mathbb{R}^m_+\setminus \{0\}$ such that the inequality~\eqref{equ-2} is fulfilled for every $y\in K$.

The proof is complete. $\hfill\Box$

\medskip
Then, for $\varepsilon$-efficient solutions of~{\rm (VP)}, the following result holds.

\begin{Theorem}\label{Proposition-2}  Let $\varepsilon\in\mathbb{R}^m_+$ and $\bar x\in K$. If $\bar x$ is an $\varepsilon$-efficient solution of~{\rm (VP)}, then there exists a vector $\lambda=(\lambda_1,\ldots,\lambda_m)\in\mathbb{R}^m_+\setminus \{0\}$ such that~\eqref{equ-2} holds. Conversely, if~\eqref{equ-2} is fulfilled for some $\lambda=(\lambda_1,\ldots,\lambda_m)\in {\rm int}\,\mathbb{R}^m_+$, then $\bar x$ is an $\varepsilon$-efficient solution of~{\rm (VP)}.
\end{Theorem}
\proof\ Since $E_{\varepsilon}\subset E_{\varepsilon}^{w}$, the first assertion follows from Theorem \ref{Proposition-1}. To prove the second assertion, suppose to the contrary that there exists a vector $\lambda=(\lambda_1,\ldots,\lambda_m)$ with $\lambda_i>0$ for all $i\in I$ such that the condition~\eqref{equ-2} is satisfied for all $y\in K$, but $\bar x\notin E_\varepsilon$. Then, by Remark~\ref{Remark 1}, there exists a vector $v\in \big[f(K)-(f(\bar x)-\varepsilon)\big]$ where $v=(v_1,...,v_m)$, $v_i\leq 0$ for all~$i\in I$, and one has $v_{i_0}<0$ for some $i_0\in I$. So, there exists $y\in K$ such that $v=f(y)-(f(\bar x)-\varepsilon)$ and 
\begin{equation}\label{system_2}
	\begin{cases}
		f_i(y)-\big(f_i(\bar x)-\varepsilon_i\big) & \leq 0\ \ \forall i\in I\setminus{\{i_0\}},\\
		f_{i_0}(y)-\big(f_{i_0}(\bar x)-\varepsilon_{i_0}\big) & < 0. 
	\end{cases}
\end{equation} 
By Lemma~\ref{Lemma1_lff}, for every $i\in I$ one has
$$f_i(y)-f_i(\bar x)=\frac{b_i^T\bar x+\beta_i}{b_i^T y+\beta_i}\langle \nabla f_i(\bar x), y-\bar x\rangle.$$ 
Hence, the system \eqref{system_2} can be rewritten as follows:
\begin{equation*}
	\begin{cases}
		\dfrac{b_i^T\bar x+\beta_i}{b_i^T y+\beta_i}\langle \nabla f_i(\bar x), y-\bar x\rangle+\varepsilon_i& \leq 0\ \ \forall i\in I\setminus{\{i_0\}},
		\\
		\dfrac{b_{i_0}^T\bar x+\beta_{i_0}}{b_{i_0}^T y+\beta_{i_0}}\langle \nabla f_{i_0}(\bar x), y-\bar x\rangle+\varepsilon_{i_0}&< 0.  
	\end{cases}
\end{equation*} 
Since $b_i^T y+\beta_i>0$ for all $i\in I$, the last system is equivalent to
\begin{equation*}\label{system-2}
	\begin{cases}
		(b_i^T\bar x+\beta_i)\langle \nabla f_i(\bar x), y-\bar x\rangle+\varepsilon_i(b_i^T y+\beta_i)&\leq 0\ \ \forall i\in I\setminus{\{i_0\}},
		\\
		(b_{i_0}^T\bar x+\beta_{i_0})\langle \nabla f_{i_0}(\bar x), y-\bar x\rangle+\varepsilon_{i_0}(b_{i_0}^T y+\beta_{i_0})&< 0.  
	\end{cases} 
\end{equation*} As $\lambda_i>0$ for all $i\in I$, one has 
\begin{equation*}\label{system-3}
	\begin{cases}
		\lambda_i\big[(b_i^T\bar x+\beta_i)\langle \nabla f_i(\bar x), y-\bar x\rangle+\varepsilon_i(b_i^T y+\beta_i)\big]&\leq 0 \ \ \forall i\in I\setminus{\{i_0\}},
		\\
		\lambda_{i_0}\big[(b_{i_0}^T\bar x+\beta_{i_0})\langle \nabla f_{i_0}(\bar x), y-\bar x\rangle+\varepsilon_{i_0}(b_{i_0}^T y+\beta_{i_0})\big]&< 0.  
	\end{cases} 
\end{equation*} Therefore, summing up the last inequalities yields
\begin{equation*}\label{equ-3}
\sum_{i=1}^{m}\lambda_i[(b_i^T\bar x+\beta_i)\langle\nabla f_i(\bar x), y-\bar x\rangle+ \varepsilon_i(b_i^Ty+\beta_i)]< 0.
\end{equation*} This obviously contradicts the condition~\eqref{equ-2}.

The proof is complete. $\hfill\Box$.

\medskip
From the above results, we can easily get necessary and sufficient conditions for $\varepsilon$-weakly efficient solution of a linear vector optimization problem as follows.

\begin{Theorem}\label{awe for lvop}
Let $\varepsilon\in\mathbb{R}^m_+$ and $\bar x\in K$. Then, $\bar x$ is an $\varepsilon$-weakly efficient solution of~{\rm (LVOP)} if and only if there exists a vector $\lambda=(\lambda_1,\ldots,\lambda_m)\in\mathbb{R}^m_+\setminus \{0\}$ such that
\begin{equation}\label{conditions for awe to lvop}
	\sum_{i=1}^{m}\lambda_i\big(\langle a_i, y-\bar x\rangle+ \varepsilon_i\big)\geq 0
\end{equation}
for all $y\in K$.
\end{Theorem}

For $\varepsilon$-efficient solutions of~{\rm (LVOP)}, we have the following result.

\begin{Theorem}\label{ae for lvop}  Let $\varepsilon\in\mathbb{R}^m_+$ and $\bar x\in K$. If $\bar x$ is an $\varepsilon$-efficient solution of~{\rm (LVOP)}, then there exists a vector $\lambda=(\lambda_1,\ldots,\lambda_m)\in\mathbb{R}^m_+\setminus \{0\}$ such that~\eqref{conditions for awe to lvop} holds. Conversely, if~\eqref{conditions for awe to lvop} is fulfilled for some $\lambda=(\lambda_1,\ldots,\lambda_m)\in {\rm int}\,\mathbb{R}^m_+$, then $\bar x$ is an $\varepsilon$-efficient solution of~{\rm (LVOP)}.
\end{Theorem}

The linear vector optimization problem in the next example has no weakly efficient solution, while both approximate efficient solution set and approximate weakly efficient solution set are nonempty.

\begin{Example}\label{T_Exam1} {\rm (See Figures~\ref{Fig 1} and \ref{Fig 2})} {\rm Consider problem $({\rm VP})$ with $n=m=2$, $$K=\{x=(x_1,x_2)\in \mathbb{R}^2\,:\,x_1\geq 0\},$$ and $f(x)=(x_1, x_2)$. One has $E^w=E=\emptyset$. Let $\varepsilon=(\varepsilon_1, \varepsilon_2)\in\mathbb{R}^2_+$ be such that $\varepsilon_1>0$. Using Remark~\ref{Remark 1}, we get $E_{\varepsilon}=\{x=(x_1,x_2)\in\mathbb R^2\,:\, 0\leq x_1<\varepsilon_1\}$ and
		$$E_{\varepsilon}^w=\{x=(x_1,x_2)\in\mathbb R^2\,:\, 0\leq x_1\leq\varepsilon_1\}.$$
		\begin{figure}[h]
			\begin{center}
				\includegraphics[height=6cm,width=9cm]{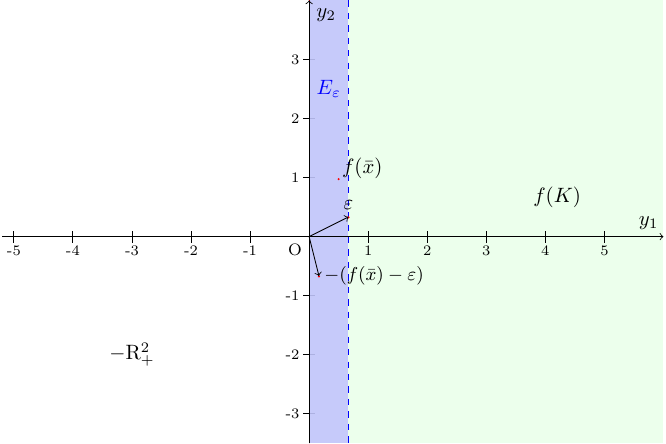}
				\caption{The set $E_{\varepsilon}$ in Example~\ref{T_Exam1}.}\label{Fig 1}
			\end{center}
		\end{figure}  
		
		\begin{figure}[h]
			\begin{center}
				\includegraphics[height=6cm,width=9cm]{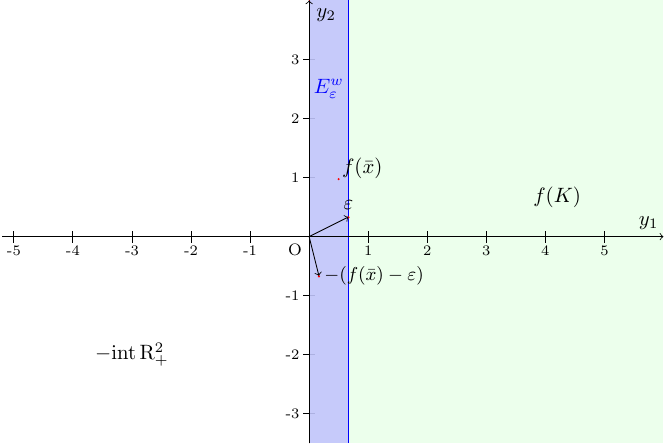}
				\caption{The set $E_{\varepsilon}^w$ in Example~\ref{T_Exam1}.}\label{Fig 2}
			\end{center}
	\end{figure} 
We can check that $\bar x=(0,0)$, $\varepsilon_1>0$ is an $\varepsilon$-weakly efficient solution by using Theorem~\ref{awe for lvop}. Indeed, choosing $\lambda=(\lambda_1,0)\in\mathbb{R}^2_+\setminus\{0\}$, $\lambda_1>0$, then the condition~\eqref{conditions for awe to lvop} is satisfied.
 Now, we assume that there exist some $\lambda=(\lambda_1,\lambda_2)\in{\rm int}\,\mathbb{R}^2_+$ such that the condition \eqref{conditions for awe to lvop} holds. This means that $$\lambda_1\big(\langle a_1, y-\bar x\rangle+ \varepsilon_1\big)+\lambda_2\big(\langle a_2, y-\bar x\rangle+ \varepsilon_2\big)\geq 0$$ for all $y\in K$.  As $a_1=(1.0)$ and $a_1=(0.1)$, one has
\begin{equation}\label{inequ_1 example}
\lambda_1\big(y_1+ \varepsilon_1\big)+\lambda_2\big(y_2+ \varepsilon_2\big)\geq 0
\end{equation} for all $y\in K$. Since $\lambda=(\lambda_1,\lambda_2)\in{\rm int}\,\mathbb{R}^2_+$  we can rewrite \eqref{inequ_1 example} equivalently as
\begin{equation}\label{inequ_2 example}
-y_2\leq \dfrac{\lambda_1\big(y_1+ \varepsilon_1\big)+\lambda_2\varepsilon_2}{\lambda_2}
\end{equation}
 for all $y\in K$.
We observer that \eqref{inequ_2 example} does not hold for all $y=(y_1,y_2)\in K, y_2<0$ and for any $\lambda\in{\rm int}\,\mathbb{R}^2_+$. This means that there is no any $\lambda=(\lambda_1,\lambda_2)\in{\rm int}\,\mathbb{R}^2_+$ such that~\eqref{ae for lvop} is satisfied. But $(0,0)$ is sill an $\varepsilon$-efficient solution by Remark~\ref{Remark 1}. Thus, Theorem~\ref{ae for lvop} cannot be used to assure that $\bar x\in E_{\varepsilon}$.}
\end{Example}

\bigskip
\noindent
\textbf{Acknowledgements} The author would like to thank Professor Nguyen Dong Yen for helpful discussions on the subject.

\end{document}